\documentclass{tac}

\usepackage{enumerate}
\usepackage{amsmath,xspace,amssymb,mathrsfs}

\usepackage{color}
\definecolor{darkgreen}{rgb}{0,0.45,0} 
\usepackage[pagebackref,colorlinks,citecolor=darkgreen,linkcolor=darkgreen]{hyperref}

\title{On the axioms for adhesive and quasiadhesive categories}

\author{Richard Garner and Stephen Lack}
\address{\\{3pt}
Computing Department Macquarie University, NSW 2109 Australia.\\
and\\{3pt}
Mathematics Department Macquarie University, NSW 2109 Australia.}

\eaddress{richard.garner@mq.edu.au, steve.lack@mq.edu.au}

\keywords{adhesive category, quasiadhesive category, pushout, exactness condition, embedding theorem}
\copyrightyear{2012}
\amsclass{18A30, 18B15}

\input xy
\xyoption{all}
\CompileMatrices

\newcommand{\rmadh}{rm-adhesive\xspace}
\newcommand{\qadh}{rm-quasiadhesive\xspace}

\newcommand{\C}{{\ensuremath{\mathscr C}}\xspace}

\newcommand{\E}{{\ensuremath{\mathscr E}}\xspace}

\newcommand{\op}{\ensuremath{{}^{\textrm{op}}}\xspace}
\newcommand{\Set}{\textnormal{\bf Set}\xspace}

\DeclareMathOperator\Sh{Sh}

  \newtheorem{proposition}{Proposition}[section]
  \newtheorem{lemma}[proposition]{Lemma}
  \newtheorem{corollary}[proposition]{Corollary}
  \newtheorem{theorem}[proposition]{Theorem}

  \theoremstyle{definition}
  \newtheorem{definition}[proposition]{Definition}
  \newtheorem{example}[proposition]{Example}

\begin{document}




\maketitle

\begin{abstract}
A category is adhesive if it has all pullbacks, all push\-outs along 
monomorphisms,
and all exactness conditions between pullbacks and pushouts along
monomorphisms which hold in a topos. This condition can be modified 
by considering only pushouts along regular monomorphisms, or by asking
only for the exactness conditions which hold in a quasitopos. We prove 
four characterization theorems dealing with adhesive categories and their
variants. 
\end{abstract}

\section{Introduction}

In the paper \cite{lex-colimits}, a general framework was given for describing
categorical structures consisting of the existence of finite limits as
well as certain types of colimits, along with exactness conditions 
stating that the limits and colimits interact in the same way as they 
do in a topos. These exactness conditions always include the fact that
the colimits in question are stable under pullback; usually, there are 
further conditions. We may sometimes speak of the colimits in question
being ``well-behaved'' when the corresponding exactness conditions
hold, it being understood that the paradigm for good behaviour is that 
which occurs in a topos. 

Examples of colimit types considered in \cite{lex-colimits} 
include coequalizers of kernel pairs,
coequalizers of equivalence relations, and finite coproducts; the
corresponding exactness conditions are then regularity, Barr-exactness,
and lextensivity. 
But for this paper, the key example from \cite{lex-colimits} is that in which 
the colimits in question are the pushouts along monomorphisms;  
the resulting categorical structure turns out to be 
that of an {\em adhesive category}, first introduced in \cite{adh}. Any 
small adhesive category admits a structure-preserving full embedding into
a topos: this was proved in \cite{adh-emb}.
Since every topos is adhesive, and every 
small adhesive category admits a full structure-preserving embedding into
a topos, we see that the exactness conditions between pushouts along 
monomorphisms and finite limits 
which hold in an adhesive category are indeed
precisely those which hold in any topos. In fact terminal objects are not 
assumed in the definition of adhesive category, and we work with pullbacks
rather than all finite limits. This does not affect the issues of limit-colimit
compatibility. 

Adhesive categories were first introduced in
order to simplify and unify various structures considered in the area of 
{\em graph transformations} \cite{FundamentalsAlgebraicGraphTransformation}, and it was clear from the start that some of the categories of interest
failed to be adhesive but still shared many of the same properties. 
There are two ways in which one could weaken the notion of adhesive 
category so as to try to include more examples. 

One possibility is
to restrict the class of colimits which are required to be well-behaved. 
We shall characterize those  categories with pullbacks and
with pushouts along regular monomorphisms,  in which the same 
exactness conditions hold as in a topos. They turn out to be precisely
those categories which were called {\em quasiadhesive} in \cite{qadh,JLS}.
The name quasiadhesive comes from the fact
that a quasiadhesive category is defined similarly to an adhesive category
except that one uses regular monomorphisms in place of monomorphisms,
just as a quasitopos is defined 
similarly to a topos except that one uses regular monomorphisms in place of 
monomorphisms. But in fact this analogy is misleading, since quasiadhesive 
categories have little to do with quasitoposes; in particular, as shown in 
\cite{JLS}, a quasitopos need not be quasiadhesive. Rather than quasiadhesive,
we shall use the name {\em \rmadh} for the categories occurring in this
characterization; see below. 

A second possibility would be to relax the requirements of good behaviour
for the colimits in question. No doubt there are many possible such 
relaxations, but we shall consider just one: requiring only 
those exactness conditions which hold in
a quasitopos. In fact we shall have nothing to say about the categories
with pullbacks, pushouts along {\em all} monomorphisms, and the exactness 
conditions holding in a quasitopos. But we shall introduce a further structure,
that of an  {\em \qadh category}, which will turn out to be that obtained 
by applying both weakenings: thus a category is \qadh just when it admits
pushouts along regular monomorphisms and pullbacks, with the same exactness properties holding between these as hold in a quasitopos.

In this paper, we prove four theorems about categorical structures involving
pushouts along some class of monomorphisms (either all of them or just
the regular ones) satisfying exactness conditions of either the stronger 
or the weaker type. 
We devote the remainder of the introduction to an explanation of the content
and significance of these theorems.

We have discussed above
notions of adhesive, \rmadh, and \qadh category in terms of exactness conditions
holding in a topos or quasitopos, but the actual definitions of these
notions are elementary, asserting properties of certain kinds of diagram
in such categories. 
A category with pullbacks and with pushouts along monomorphisms is
said to be adhesive when the pushouts along monomorphisms are 
{\em van Kampen}; the name coming from the fact that this condition is
formally similar to the statement of  
the coverings form of the van Kampen theorem 
\cite{Brown-Janelidze}. Explicitly, the condition says that for any cube 
$$\xymatrix @R1pc @C1pc { 
& C' \ar[dl]_{m'} \ar'[d][dd]_c \ar[rr]^{f'} && B' \ar[dl]^{n'} \ar[dd]^{b} \\ 
A' \ar[dd]_{a} \ar[rr]^(0.7){g'} && D' \ar[dd]^(0.7){d} \\
& C \ar'[r][rr]^(0.5){f} \ar[dl]_{m} && B \ar[dl]^{n} \\
A \ar[rr]_{g} & & D }$$
with the bottom face as the given square and with the left and back faces
as pullbacks, the top face is a pushout if and only if the front and right 
faces are pullbacks. An elegant reformulation of the van Kampen condition
was given in \cite{HeindelSobocinski}, where it was shown that a pushout
along a monomorphism is van Kampen just when it is also a (bicategorical)
pushout in the bicategory of spans. 

The ``if'' part of the van Kampen condition says that the pushout
is stable under pullback. We shall from now on use the adjective
``stable'' to mean ``stable under pullback''. It is easy to see that a 
van Kampen square is a pullback as well as a stable pushout; in general,
however, there can be stable pushouts which are pullbacks but fail to 
be van Kampen. 

We now collect together the elementary definitions of our various notions. 
\begin{definition}
A category with pullbacks is said to be: 
\begin{enumerate}[(i)]
\item {\em adhesive} if it has pushouts along monomorphisms and these are van Kampen;
\item {\em \rmadh} if it has pushouts along regular monomorphisms and these are van Kampen;
\item {\em \qadh} if it has pushouts along regular monomorphisms and these are stable and are pullbacks.
\end{enumerate}
\end{definition}

It might appear that we should have considered a fourth possibility,
namely the categories which  have pushouts along monomorphisms
that are stable and are pullbacks, but need not be van Kampen.
Our first theorem states that in this case the pushouts along monomorphisms must in fact be van Kampen, and so this apparently different 
notion coincides with that of adhesive category.
It will be proved as part of Theorem~\ref{thmA} below. 


\subsection*{Theorem~A}
A category with pullbacks is adhesive if and only if it has pushouts along
monomorphisms, and these pushouts are stable and are pullbacks.

\medskip

In fact this characterization was observed in \cite{lex-colimits}, with an 
indirect proof based on the fact that the proof of the  embedding theorem for 
adhesive categories given in \cite{adh-emb} used only the elements of
the new characterization, not the van Kampen condition. The proof given 
below, however, is direct. If we take the equivalent condition given in the
theorem as the definition of adhesive category, then the fact that every 
topos is adhesive \cite{topos-adh,adh-emb} becomes entirely standard. 

Our second main result, proved as Corollary~\ref{thmB} below, 
describes the relationship between \rmadh and \qadh categories:

\subsection*{Theorem~B}
A category is \rmadh just when it is \qadh and regular subobjects are 
closed under binary union.

\medskip

This is a sort of refinement of the result in \cite{JLS} which states that a 
quasitopos is \rmadh (there called quasiadhesive) if and only if regular
subobjects are closed under binary union. As observed in 
\cite[Corollary~20]{JLS}, the category of sets equipped with a 
binary relation is a quasitopos, and so \qadh, but fails to be \rmadh.

The remaining two theorems are embedding theorems, and complete 
the characterizations of \rmadh and \qadh categories as those with 
pushouts along regular monomorphisms satisfying the exactness 
conditions which hold in a topos or a quasitopos.

\subsection*{Theorem~C}
Let \C be a small category with all pullbacks and with pushouts along
regular monomorphisms. Then \C is \rmadh just when it has
a full embedding, preserving the given structure, into a topos. 

\medskip

It follows from the theorem that the condition of being \rmadh can be 
seen as an exactness condition in the sense of \cite{lex-colimits}. 
The  theorem also clinches the argument that the name quasiadhesive 
does not really belong with the condition that all pushouts along regular
monomorphisms be van Kampen; this condition is related to toposes
rather than quasitoposes. 

\subsection*{Theorem~D}
Let \C be a small category with all pullbacks and with all pushouts 
along regular monomorphisms. Then \C is \qadh just when it has a full 
embedding, preserving the given structure, into a quasitopos.

\medskip

One direction of each theorem is immediate, since every topos is
adhesive and every quasitopos is \qadh. The converse implications 
are proved in Theorems~\ref{thmC}
and~\ref{thmD} respectively; the theorems themselves are stated as
Corollaries~\ref{corC} and~\ref{corD}.

We have found it convenient to give the name {\em adhesive morphism} 
to a morphism $m$ (necessarily a regular monomorphism) with the property
that pushouts along any pullback of $m$ are stable and are pullbacks. 
One aspect of the convenience is that we can then deal simultaneously
with the context where all monomorphisms are adhesive and that where
all regular monomorphisms are adhesive. Other authors have defined
adhesiveness with respect to an abstract class of morphisms---see
\cite{FundamentalsAlgebraicGraphTransformation} and the references 
therein---but this is subtly different. Rather than specifying a class 
of morphisms and asking that pushouts along them be suitably well 
behaved, we are instead considering the class of all morphisms 
along which pushouts are suitably well behaved. 

We introduce and study these adhesive morphisms in 
Section~\ref{sect:adh}. Then in Section~\ref{sect:VK} we look at 
when pushouts along adhesive morphisms are van Kampen squares, and 
prove the first two theorems described above. Finally in Section~\ref{sect:emb}
we look at structure preserving embeddings of \rmadh and \qadh 
categories, and prove the last two theorems.

\subsection*{Acknowledgements} 
The second-named author is grateful to Peter Johnstone and Pawe\l\
Soboci\'nski for helpful conversations in relation to the content of the 
paper. 

This research was supported under the Australian Research
Council's {\em Discovery Projects} funding scheme, 
project numbers DP110102360 (Garner) and DP1094883 (Lack). 


\section{Adhesive morphisms}\label{sect:adh}

Let \C be a category with pullbacks. 
Say that a morphism $m$ is {\em pre-adhesive}, if pushouts along $m$
exist, are stable, and are pullbacks. Say that $m$ is {\em adhesive} if 
all of its pullbacks are pre-adhesive.

\begin{proposition}
All isomorphisms are adhesive. 
All adhesive morphisms are pre-adhesive. 
  All pre-adhesive morphisms are regular mono\-morphisms. 
\end{proposition}

\proof 
Any isomorphism is obviously pre-adhesive,
but the pullback of an isomorphism is still an isomorphism, and so 
isomorphisms are adhesive. The second statement is a triviality.
As for the third, if $m$ is pre-adhesive, then it is the pullback of its pushout 
along itself. This says that it is the equalizer of its cokernel pair, and
so is a regular monomorphism. 
\endproof

We shall repeatedly use the standard results about pasting and cancellation
of pullbacks and pushouts: if in a diagram
$$\xymatrix{
\ar[r] \ar[d] & \ar[r] \ar[d] & \ar[d] \\
\ar[r]_r & \ar[r] & }$$
the right square is a pullback, then the left square is a pullback just when 
the composite rectangle is one; dually, if the left square is a pushout, then
the right square is a pushout just when the composite rectangle is one. 
It is well known that there is another cancellation result for pullbacks:
 if the left square and the composite rectangle are pullbacks,
and the arrow $r$ is a stable regular epimorphism, then the right square
is a pullback. This is proved in 
\cite[Lemma~4.6]{CJKP-LocalizationStabilization} under the assumption
that $r$ is an effective descent morphism, but also follows from the 
special case $a_1=a_2=r$ of the 
following closely related result, which we shall use repeatedly.


\begin{lemma}\label{lemma:cancellation}
Consider the following diagrams
$$\xymatrix{ 
A_0 \ar[r]^{b_2} \ar[d]_{b_1} & A_2 \ar[d]^{a_2} \\
A_1 \ar[r]_{a_1} & A }
\xymatrix{
A'_i \ar[d]_{p_i} \ar[r]^{a'_i} & A' \ar[d]^{p} \ar[r]^{f'} & B' \ar[d]^{q} \\
A_i \ar[r]_{a_i} & A \ar[r]_{f} & B }$$
 in a category \C with pullbacks. Suppose that 
the square on the left is a stable pushout, and the 
the two diagrams on the right ($i=1,2$) have left square a 
pullback and composite a pullback. Then the square on the 
right is also a pullback.
\end{lemma}

\proof
Consider the diagrams ($i=1,2$)  
$$\xymatrix{
A'_0 \ar[r]^{b'_i} \ar@{=}[d] & A'_i \ar[r]^{a'_i} \ar@{=}[d] & 
A' \ar[r]^{f'} \ar[d]_{p'} & B' \ar@{=}[d] \\
A'_0 \ar[r]^{b'_i} \ar[d]_{p_0} & 
A'_i \ar[r]^{a''_i} \ar[d]_{p_i} & 
A'' \ar[r]^{f''} \ar[d]_{q'} & 
B' \ar[d]^{q} \\
A_0 \ar[r]_{b_i} & A_i \ar[r]_{a_i} & A \ar[r]_{f} & B }$$
in which the squares along the bottom are all pullbacks, and 
$p'$ is the unique morphism satisfying $f''p'=f'$ and $q'p'=p$. 
We need to show that $p'$ is invertible. 

By the standard properties of pullbacks, the left and middle 
squares on the top row are pullbacks; thus in the diagram 
$$\xymatrix{
A'_0 \ar[r]^{b'_2} \ar[d]_{b'_1} & A'_2 \ar[d]_{a'_2} \ar[ddr]^{a''_2} \\
A'_1 \ar[r]^{a'_1} \ar[drr]_{a''_1} & A' \ar[dr]^(0.2){p'} \\
&& A'' }$$
both the interior and the exterior square are pullbacks of a stable 
pushout, hence are pushouts. Since $p'$ is the canonical 
comparison between two pushouts, it is invertible. \endproof

The class of adhesive morphisms is of course stable under pullback, 
since it consists of precisely the stably-pre-adhesive morphisms. We 
also have:

\begin{proposition}
The class of adhesive morphisms is closed under composition and stable 
under pushout.
\end{proposition}

\proof
Closure under composition follows immediately from the pasting 
properties of pullbacks and pushouts.

As for stability, suppose that $m\colon C\to A$ is adhesive, and consider a
pushout
$$\xymatrix{
C \ar[d]_{m} \ar[r]^{f}  & B \ar[d]^{n} \\
A \ar[r]_{g} & D }$$ 
which is also a pullback, since $m$ is adhesive. 
First we show that $n$ is mono. Pull back this pushout along $n$ to get a
cube
$$\xymatrix @R1pc @C1pc { 
& C \ar@{=}[dl] \ar@{=}'[d][dd] \ar[rr] && B' \ar[dl]^{n_1} \ar[dd]^{n_2} \\ 
C \ar[dd]_{m} \ar[rr]^(0.7){f}   && B \ar[dd]^(0.7){n} \\
& C \ar'[r][rr]^(0.5){f} \ar[dl]_{m} && B \ar[dl]^{n} \\
A \ar[rr]_{g} & & D }$$
in which the top face, like the bottom, is a pushout, and the four
vertical faces are pullbacks; and in which  
$n_1$ and $n_2$ are the kernel pair of $n$.
Since $n_1$ is a pushout of an isomorphism it is an isomorphism. 
This implies that the kernel pair of $n$ is trivial, and so that $n$ is a 
monomorphism; thus any pushout
of an adhesive morphism is a monomorphism. 

Next we show that the monomorphism $n$ is pre-adhesive. Given any 
$h\colon B\to E$ we can push out $m$ along $hf$, and so by the cancellativity
property of pushouts obtain pushout squares
$$\xymatrix{
C \ar[d]_{m} \ar[r]^{f} & B \ar[d]_{n} \ar[r]^{h} & E \ar[d]^{p} \\
A \ar[r]_{g} & D \ar[r]_{k} & F }$$
where the pushout on the right is stable, since the left and composite 
pushouts are so. We must show that the square on the right is 
also a pullback.  Since the square on the left is a stable pushout, 
we may use Lemma~\ref{lemma:cancellation}. We know that the composite
rectangle is a pullback, so it will suffice to show that in the diagram 
below
$$\xymatrix{
B \ar@{=}[r] \ar@{=}[d] & B \ar[d]_{n} \ar[r]^{h} & E \ar[d]^{p} \\
B \ar[r]_{n} & D \ar[r]_{k} & F }$$
the left square and the composite rectangle are pullbacks. The left square
is indeed a pullback because $n$ is monic; the composite is a pullback
because $p$ is is monic.

Finally, to see that $n$ is adhesive, let $n'\colon B'\to D'$ be a pullback of $n$ along 
$d\colon D'\to D$. Then we have a cube
$$\xymatrix @R1pc @C1pc { 
& C' \ar[dl]_{m'} \ar'[d][dd]^c \ar[rr]^{f'} && B' \ar[dl]^{n'} \ar[dd]^{b} \\ 
A' \ar[dd]_{a} \ar[rr]^(0.7){g'} && D' \ar[dd]^(0.7){d} \\
& C \ar'[r][rr]^(0.5){f} \ar[dl]_{m} && B \ar[dl]^{n} \\
A \ar[rr]_{g} & & D }$$
in which the bottom face is a stable pushout since $m$ is pre-adhesive,
so that the top face is a pushout. But $m'$ is adhesive since $m$ is
adhesive, and so $n'$ is also pre-adhesive. This now proves
that $n$ is adhesive.
\endproof

Recall that a union of two subobjects is said to be {\em effective} if it 
can be constructed as the pushout over the intersection of the subobjects.
This will be the case provided that the pushout exists, and the  induced map 
out of the pushout is a monomorphism. 

\begin{proposition}\label{prop:effective}
Binary unions are effective provided that at least one of the subobjects
is adhesive.
\end{proposition}

\proof
  Consider a subobject $h\colon A\to E$ and an adhesive subobject
$p\colon B\to E$, and form their intersection $C$, and their pushout $D$
over $C$, as in the diagram 
$$\xymatrix{
C \ar[r]^{f} \ar[d]_{m} & B \ar[d]_{n} \ar[ddr]^{p} \\
A \ar[r]^{g} \ar[drr]_{h} & D \ar[dr]^(0.3){x} \\
&& E }$$
in which $x\colon D\to E$ is the induced map; we must show
that $x$ is a monomorphism. Observe that 
$g$ and $f$ are monomorphisms since $h$ is one, while $n$ and 
$m$ are adhesive since $p$ is so.

Pulling back the external part of the diagram along $h\colon A\to E$ gives
the external part of the diagram
$$\xymatrix{
C \ar@{=}[r] \ar[d]_{m} & C \ar[d]_{m} \ar[ddr]^{m} \\
A \ar@{=}[r]  \ar@{=}[drr] & A \ar@[=][dr] \\
&& A}$$
but pulling back the pushout in the interior of the original diagram
must give a pushout, which must then be the square in the interior of
the diagram above. 
 It follows that we have a pullback
$$\xymatrix{
A \ar[r]^{g} \ar@{=}[d] & D \ar[d]^{x} \\
A \ar[r]_{h} & E. }$$

Similarly, pulling back the original diagram along $p\colon B\to E$ gives the diagram
$$\xymatrix{
C \ar[r]^{f} \ar@{=}[d] & B \ar@{=}[d] \ar@{=}[ddr] \\
C \ar[r]^{f} \ar[drr]_{f} & B \ar@{=}[dr] \\
&& B }$$
and we have a pullback
$$\xymatrix{
B \ar[r]^{n} \ar@{=}[d] & D \ar[d]^{x} \\
B \ar[r]_{p} & E. }$$

We now use Lemma~\ref{lemma:cancellation} to prove that 
$$\xymatrix{
D \ar@{=}[r] \ar@{=}[d] & D \ar[d]^{x} \\
D \ar[r]_{x} & E }$$
is a pullback, and so that $x$ is a monomorphism, using the stable pushout
$$\xymatrix{
C \ar[r]^{f} \ar[d]_{m} & B \ar[d]^{n} \\
A \ar[r]_{g} & D.}$$
This is done by observing that in the rectangles
$$\xymatrix{
B \ar[r]^{n} \ar@{=}[d] & D \ar@{=}[r] \ar@{=}[d] & D \ar[d]^{x} \\
B \ar[r]_{n} & D \ar[r]_{x} & E } \quad
\xymatrix{
A \ar[r]^{g} \ar@{=}[d] & D \ar@{=}[r] \ar@{=}[d] & D \ar[d]^{x} \\
A \ar[r]_{g} & D \ar[r]_{x} & E }
$$
the left square and the composite are pullbacks, hence the (common) 
right square is one. 
\endproof

The following lemma will be used several times in later sections of the paper;
the argument in the proof  goes back to \cite{topos-adh}, where it was
used in the context of adhesive categories; for a general framework which 
explains the origin of the condition itself, see \cite[Section~6]{lex-colimits}.

\begin{lemma}\label{lemma:basic}
Suppose that adhesive subobjects are closed under binary union, and that
all split monomorphisms are adhesive. 
Let $m\colon C\to A$ be an adhesive morphism and $f\colon C\to B$ any morphism.
Construct the diagram
$$\xymatrix{
C \ar[r]^{\gamma} \ar[d]_m & C_2 \ar[d]_{m_2} \ar@<1ex>[r]^{f_1}\ar@<-1ex>[r]_{f_2}
& C \ar[r]^{f} \ar[d]_m & B \ar[d]^n \\
A \ar[r]_{\delta} & A_2 \ar@<1ex>[r]^{g_1} \ar@<-1ex>[r]_{g_2} & A \ar[r]_g & D }$$
in  which the right hand square is a pushout, $(f_1,f_2)$ and $(g_1,g_2)$ are the
kernel pairs of $f$ and $g$, with $m_2$ the induced map, 
while $\gamma$ and $\delta$ are the diagonal maps.
Then the left square, the right square, and the two central squares
are all both pushouts and pullbacks.
\end{lemma}

\proof
Since $m$ is adhesive, the right hand square is a stable pushout. We may pull it back along $g$ so as to give one of the central squares
$$\xymatrix{
C_2 \ar[r]^{f_1} \ar[d]_{m_2} & C \ar[d]^{m} \\
A_2 \ar[r]_{g_1} & A }$$
which is therefore a pushout and a pullback; 
the case of the other central square
is similar. Furthermore $m_2$ is adhesive, since it is a pullback of $m$. 
Form the diagram
$$\xymatrix{
C \ar[r]^\gamma \ar[d]_m & C_2 \ar[r]^{f_1} \ar[d]^j & C \ar[d]^m \\
A \ar[r]^i \ar[dr]_\delta & E\ar[r]^{h_1} \ar[d]^k & A \ar[d]^1 \\
& A_2 \ar[r]_{g_1} & A }$$
in which the top left square is a pushout, $h_1$ is the unique morphism satisfying
$h_1i=1$ and $h_1j=mf_1$, and 
$k$ is the unique morphism satisfying $ki=\delta$ and $kj=m_2$. 
By Proposition~\ref{prop:effective}, the morphism $k\colon E\to A_2$ is the 
union of $\delta\colon A\to A_2$ and the adhesive subobject 
$m_2\colon C_2\to A_2$; since $\delta$ is split, it is adhesive by one 
of our hypotheses on the category, thus the union $k$ of $\delta$ 
and $m_2$ is adhesive by the other hypothesis. 
We are to show that $k$ is invertible.

Now the top left square and the composite of the upper squares are 
both pushouts, so the top right square is also a pushout, by the cancellativity
properties of pushouts. The composite of the two squares on the right is the pushout constructed at the beginning of the proof,
so finally the lower right square is a pushout by the cancellativity property of 
pushouts once again.

Since $k$ is adhesive, the lower right square is a pullback.
Thus $k$ is invertible, and our square is indeed
a pushout. 
\endproof

\section{The van Kampen condition}\label{sect:VK}

We recalled in the introduction the notion of a van Kampen square, which
appears in the definition of adhesive and \rmadh category. In this 
section we give various results showing how this condition can be
reformulated. 

\begin{proposition}\label{prop:adhesive}
  Suppose that adhesive subobjects are closed under binary union
and that split subobjects are adhesive. Then
pushouts along adhesive morphisms are van Kampen.
\end{proposition}

\proof
Consider a cube 
   $$\xymatrix @R1pc @C1pc {
 & C' \ar[rr]^{f'} \ar[dl]_{m'} \ar'[d][dd]^c && B' \ar[dl]_{n'} \ar[dd]^b \\
 A' \ar[rr]^(0.7){g'} \ar[dd]_a && D' \ar[dd]^(0.7)d \\
 & C \ar'[r][rr]^(0.5)f \ar[dl]_m && B \ar[dl]^n \\
 A \ar[rr]_g && D }$$
in which $m$ is adhesive, the bottom face is a pushout, and the left 
and back faces are pullbacks. We are to show that the top face is a 
pushout if and only if the front and right faces are pullbacks. 

Since $m$ is adhesive, the pushout on the bottom face is stable, thus 
if the front and right faces are pullbacks, then the top face is also a 
pushout. 

Suppose conversely that the top face is a pushout, noting also that 
$m'$ is a pullback of the adhesive morphism $m$ so is itself adhesive,
and so that $n'$ and $n$ are also adhesive. 
We are to show that the front and right faces are pullbacks.
We do this using various applications of Lemma~\ref{lemma:cancellation}.

To see that the right face is a pullback, apply Lemma~\ref{lemma:cancellation} with the top face as stable pushout. 
We know that the top face, the left face and the bottom face
are pullbacks; thus also the composite of top and right faces are
pullbacks. Thus it remains to show that the composite and the 
left square in the diagram on the
left 
$$\xymatrix{
B' \ar@{=}[r] \ar@{=}[d] & B' \ar[r]^{b} \ar[d]_{n'} & B \ar[d]^{n} \\
B' \ar[r]_{n'} & D' \ar[r]_{d} & D }\quad
\xymatrix{
B' \ar[r]^{b} \ar@{=}[d] & B \ar@{=}[r] \ar@{=}[d] & B \ar[d]^{n} \\
B' \ar[r]_{b} & B \ar[r]_{n} & D 
}$$
are pullbacks. Of these, the left square
is a pullback since $n'$
is monic, while the composite is equally the composite of the
diagram on the right, in which both squares are pullbacks
(because $n$ is monic). 

This proves that the right face is a pullback; now consider the
front face  
$$\xymatrix{
A' \ar[d]_{g'} \ar[r]^{a} & A \ar[d]^{g} \\
D' \ar[r]_{d} & D }$$
(here reflected about the diagonal).
We shall show that it is a pullback using Lemma~\ref{lemma:cancellation} with the stable pushout 
$$\xymatrix{
C' \ar[r]^{m'} \ar[d]_{f'} & A' \ar[d]^{g'} \\
B' \ar[r]_{n'} & D'. }$$
We know that the back and bottom faces of the cube are pullbacks,
and so the composite of the top and front faces is a pullback.
Thus it will suffice to show that if, in the diagram below on the left,
$$
\xymatrix{
A'_2 \ar[r]^{g'_1} \ar[d]_{g'_2} & A' \ar[r]^{a} \ar[d]_{g'} & A \ar[d]^{g} \\
A' \ar[r]_{g'} & D' \ar[r]_{d} & D }\quad
\xymatrix{
A'_2 \ar[r]^{a_2} \ar[d]_{g'_2} & A_2 \ar[r]^{g_1} \ar[d]_{g_2} & 
A \ar[d]^{g} \\
A' \ar[r]_{a} & A \ar[r]_{g} & D }
$$
the left square is a pullback, then so is the composite. To see that this is so,
consider the diagram on the right in which the right square is a pullback,
and $a_2$ is the unique morphism making the square on the left commute, and the composite rectangle equal the composite rectangle on the left.
It will now suffice to show that the left square of the right diagram 
(here reflected in the diagonal) 
$$\xymatrix{
A'_2 \ar[d]_{a_2} \ar[r]^{g'_2} & A' \ar[d]^{a} \\
A_2 \ar[r]_{g_2} & A }$$
is a pullback. We do this using Lemma~\ref{lemma:cancellation} once
again, this time using the  stable pushout 
$$\xymatrix{
C \ar[r]^{\gamma} \ar[d]_{m} & C_2 \ar[d]^{m_2} \\
A \ar[r]_{\delta} & A_2}$$
constructed in  Lemma~\ref{lemma:basic} (it is 
here that we need the assumption that adhesive subobjects are closed
under binary union, and contain the split subobjects).

First of all, in the diagram 
$$
\xymatrix{
A' \ar[r]^{\delta'} \ar[d]_{a} & A'_2 \ar[d]_{a_2} \ar[r]^{g'_2} & A' \ar[d]^{a} \\
A \ar[r]_{\delta} & A_2 \ar[r]_{g_2} & A,}$$
the horizontal composites are identities, and so the composite
is a pullback, while the left square is equally the right face of the 
cube
\begin{equation}
  \label{eq:1}
\vcenter{\hbox{
\xymatrix @R1pc @C1pc {
 & C' \ar[rr]^{m'} \ar[dl]_{\gamma'} \ar'[d][dd]^c && A' \ar[dl]_{\delta'} \ar[dd]^a \\
 C'_2 \ar[rr]^(0.7){m'_2} \ar[dd]_{c_2} && A'_2 \ar[dd]_(0.7){a_2} \\
 & C \ar'[r][rr]^(0.5){m} \ar[dl]_{\gamma} && A \ar[dl]^{\delta} \\
 C_2 \ar[rr]_{m_2} && A_2 }  
}}
\end{equation}
and so is a pullback by the earlier part of the proof, since in this cube
the top and bottom faces are pushouts along adhesive morphisms,
and the left and back faces are pullbacks. Since $m$ is also 
adhesive, by the earlier part of the proof once again we may deduce 
that the front square is also a pullback. 

Thus it remains to show that the composite and the
left square in the diagram on the left
$$
\xymatrix{
C'_2 \ar[r]^{m'_2} \ar[d]_{c_2} & A'_2 \ar[r]^{g'_2} \ar[d]_{a_2} & A' \ar[d]^{a} \\
C_2 \ar[r]_{m_2} & A_2 \ar[r]_{g_2} & A }\quad
\xymatrix{
C'_2 \ar[r]^{f'_2} \ar[d]_{c_2} & C' \ar[r]^{m'} \ar[d]_{c} & A' \ar[d]^{a} \\
C_2 \ar[r]_{f_2} & C \ar[r]_{m} & A }
$$
are pullbacks. Of these,  the left square is equally the 
front face of the cube~\eqref{eq:1}, which we already saw to be a pullback. 
On the other hand, the  composite is equally the composite of the 
diagram on the right, in which both squares are pullbacks.
\endproof

We are now ready to give our first main theorem, providing a new 
characterization of adhesive categories. 

\begin{theorem}\label{thmA}
For a category \C with pullbacks, the following conditions are equivalent:
\begin{enumerate}[(i)]
\item \C is adhesive;
\item all monomorphisms are adhesive;
\item \C has pushouts along monomorphisms, and these pushouts are
stable and are pullbacks.
\end{enumerate}
\end{theorem}

\proof
Here condition ($iii$) says that all monomorphisms are pre-adhesive.
Clearly ($i$) implies ($ii$) and ($ii$) implies ($iii$). If all monomorphisms
are pre-adhesive, then since the pullback of a monomorphism is a
monomorphism, all monomorphisms are in fact adhesive; this shows
that ($iii$) implies ($ii$). Finally if all monomorphisms are adhesive,
then certainly the adhesive subobjects are closed under binary union, so 
all pushouts along monomorphisms are van Kampen by Proposition~\ref{prop:adhesive}. Thus ($ii$) implies ($i$).
\endproof

We recall from the introduction that we define a category to be \rmadh 
if it has all pullbacks
and all pushouts along regular monomorphisms, and these pushouts
are van Kampen; while if instead these pushouts are only assumed to be stable
and to be pullbacks, then the category is called \qadh. Thus a category with pullbacks is \qadh if and only if all regular monomorphisms are adhesive. We also reiterate the 
warning that the notion here being called \rmadh was formerly known as
quasiadhesive. 

\begin{theorem}
  For a category \C with pullbacks, the following conditions are equivalent:
  \begin{enumerate}[(i)]
  \item \C is \rmadh;
    \item all regular monomorphisms are adhesive, and regular subobjects
      are closed under binary union;
      \item \C has pushouts along regular monomorphisms, these pushouts are
stable and are pullbacks, and regular subobjects are closed under binary union.
  \end{enumerate}
\end{theorem}

\proof
The downward implications are all straightforward, once we note the fact,
proved in \cite{JLS}, that in an \rmadh category, the union of two regular
subobjects is again regular. The upward implications are proved just
as in the proof of the previous theorem.
\endproof

In particular, we have the following result, which refines the fact \cite{JLS}
that a quasitopos is \rmadh if and only if the union of two regular 
subobjects is again regular. 

\begin{corollary}\label{thmB}
A category is \rmadh just when it is \qadh and regular subobjects are 
closed under binary union.
\end{corollary}

Any locally presentable locally cartesian closed category of course
has all colimits and these are all stable under pullback. This is not enough
to imply that the category is \qadh, as the following
example shows.  The example was originally given in \cite{BorceuxPedicchio}, 
attributed to Ad\'amek and Rosick\'y, as an example of a locally 
presentable and locally cartesian closed category which is not a quasitopos.

\begin{example}
Let \E be the category of sets equipped with a reflexive binary relation
with the property that 
\begin{align*}
  x_1Rx_2Rx_1 &\Rightarrow x_1=x_2 \\
x_1Rx_2Rx_3Rx_1 &\Rightarrow x_1=x_2=x_3 \\
\vdots & ~~~\vdots \\
x_1Rx_2R\ldots Rx_nRx_1 &\Rightarrow x_1=x_2=\ldots=x_n
\end{align*}
for all $n\ge 2$. As observed in \cite[Section~4]{BorceuxPedicchio},
this is locally cartesian closed and locally presentable, but the 
regular monomorphism given by the inclusion $(\{0,2\},R)\to(\{0,1,2\},S)$
is not classified, where $R$ consists only of the diagonal, and where $S$
consists of the diagonal as well as $0S1$ and $1S2$. In fact the pushout
of this inclusion along the unique map from $(\{0,2\},R)$ to the terminal
object is clearly just the terminal object, and so this pushout is not 
also a pullback. Thus \E is not \qadh.
\end{example}

\section{Embeddings}\label{sect:emb}

In this section we prove various embedding theorems. In each case,
the category \C will be fully embedded into a full subcategory \E of
the presheaf category $[\C\op,\Set]$, via the corresponding restricted 
Yoneda functor; the category \E will be a topos or a quasitopos, as appropriate; and the embedding will preserve the relevant colimits. 
The Yoneda functor $\C\to[\C\op,\Set]$, of course,
preserves limits but not colimits; we fix this by  choosing \E appropriately.
The preservation 
of colimits will occur provided that all presheaves $F\colon\C\op\to\Set$
lying in \E send the relevant colimits in \C to limits in \Set. In order to 
guarantee that \E is a topos, we take it to be the category of sheaves
for a topology $j$; we choose $j$ in such a way that the representables are 
sheaves, and all sheaves send the relevant 
colimits in \C to limits in \Set. If instead we only need \E to be a quasitopos,
then we need only take \E to be the category of $j$-sheaves which are 
also $k$-separated, where $j$  and $k$ are topologies, with $k$ containing
$j$: such quasitoposes are called Grothendieck quasitoposes, and were
studied in \cite{BorceuxPedicchio}. In each case $j$ will be a 
cd-topology, in the sense of \cite{voevodsky:cd-topology}. 

Suppose that \C is a small category with pullbacks.  Since pushouts along
adhesive morphisms are stable, there is 
a topology $j$ on \C in which the covering families are those of the form
$\{g,n\}$ for some pushout
$$\xymatrix{
C \ar[r]^f \ar[d]_m & B \ar[d]^n \\ A \ar[r]_g & D }$$
with $m$ adhesive. 

The sheaves for the topology are those functors $F\colon \C\op\to\Set$ for which,
given a diagram of the above form, the induced diagram
\begin{equation}\label{sheaf}
\xymatrix{
& FC & FB \ar[l]_{Ff} \\
FA_2 & FA \ar[u]^{Fm} \ar@<1ex>[l]^{Fg_1} \ar@<-1ex>[l]_{Fg_2} & 
FD \ar[l]_{Fg} \ar[u]_{Fn} }
\end{equation}
in which $(g_1,g_2)$ is the kernel pair of $g$, exhibits $FD$ as the limit
in \Set of the remainder of the diagram.

First suppose that split monomorphisms are adhesive,
and adhesive subobjects are closed under binary union. Then by
Lemma~\ref{lemma:basic}, the diagram
$$\xymatrix{
C \ar[r]^{\gamma} \ar[d]_{m} & C_2 \ar[d]^{m_2} \\
A \ar[r]_{\delta} & A_2}$$
is a pushout along the adhesive morphism $m$, and so $\{\delta,m_2\}$
is a covering family. Since both $\delta$ and $m_2$ are monomorphisms,
any sheaf $F$ sends this pushout to a pullback in \Set, and so in 
particular $F\delta$ and $Fm_2$ are jointly monic. As explained
in \cite{adh-emb}, this now
allows the sheaf condition given in terms of the
diagram \eqref{sheaf} above to be simplified to the requirement that
the square in that diagram be a pullback; the argument is sufficiently
short that it is worth repeating here. Suppose that $x\in FA$ and $y\in FB$
are given, with $Fm.x=Ff.y$; we must show that also $Fg_1.x=Fg_2.x$, so
that by the sheaf condition there is a unique $z\in FD$ with 
$Fg.z=x$ and $Fn.z=y$. To do this, observe that 
\begin{align*}
F\delta.Fg_1.x &=F(g_1\delta).x=x=F(g_2\delta).x=F\delta.Fg_2.x \\
\intertext{and}
Fm_2.Fg_1.x &=Ff_1.Fm.x=Ff_1.Ff.y=Ff_2.Ff.y \\
&=Ff_2.Fm.x=Fm_2.Fg_2.x
\end{align*}
thus since $F\delta$ and $Fm_2$ are jointly monic, $Fg_1.x=Fg_2.x$
as required. 

The sheaves are therefore precisely those presheaves which send push\-outs
along adhesive morphisms to pullbacks. Then the restricted 
Yoneda embedding $\C\to\Sh(\C)$ preserves such pushouts as 
well as any existing limits. This proves:

\begin{theorem}
  For any small category \C with pullbacks, in which split monomorphisms
are adhesive, and adhesive subobjects are closed under binary union,
there is a full embedding of \C into a topos, and this embedding preserves
pushouts along adhesive morphisms as well as all existing limits.
\end{theorem}

Thus we recover the result of \cite{adh-emb} that any small adhesive
category has a full structure-preserving embedding in a topos. But 
we also obtain:

\begin{theorem}\label{thmC}
Any small \rmadh category can be fully embedded in a topos,
via an embedding which preserves pushouts along regular monomorphisms
as well as all existing limits.
\end{theorem}

Conversely, it is well-known that every topos satisfies the 
conditions in the definition of \rmadh category, hence so will any
full subcategory closed under the relevant limits and colimits. Thus
we have:

\begin{corollary}\label{corC}
Let \C be a small category with all pullbacks and with pushouts along
regular monomorphisms. Then \C is \rmadh just when it has
a full embedding, preserving the given structure, into a topos. 
\end{corollary}

We now turn to the embedding theorem for \qadh categories. First, however,
we need to develop the theory of such categories a little.  
The next result is essentially contained in \cite{JLS}, although the 
context there was slightly different. 

\begin{proposition}
If $m\colon A\to X$ is the union of regular subobjects $m_1\colon A_1\to X$ and
$m_2\colon A_2\to X$ in an \qadh category, then $m$ has a stable (epi,regular mono) factorization.
\end{proposition}

\proof
We construct the union as the pushout square in 
$$\xymatrix{
A_0 \ar[r]^{m'_1} \ar[d]_{m'_2} & A_2 \ar[d]_{n_2} \ar[ddr]^{m_2} \\
A_1 \ar[r]^{n_1} \ar[drr]_{m_1} & A \ar[dr]^(0.2){m} \\
&& X }$$
where $A_0$ is the intersection of $A_1$ and $A_2$. 
First let $i,j\colon X\rightrightarrows X_1$ be the cokernel pair of $m_1$, and $e_1\colon X_1\to X$
the codiagonal. We can pull all these maps back along $m_2$ to get
the diagram
$$\xymatrix{
A_0 \ar[r]^{m'_1} \ar[d]_{m'_2} & A_2 \ar[d]_{m_2} \ar@<1ex>[r]^{i_2} 
\ar@<-1ex>[r]_{j_2} & X_2 \ar[r]^{e_2} \ar[d]_{\ell} & A_2 \ar[d]^{m_2} \\
A_1 \ar[r]^{m_1} & X \ar@<1ex>[r]^{i} \ar@<-1ex>[r]_{j} & X_1 \ar[r]^{e_1} &
X }$$
in which $\ell$ is a pullback of the regular monomorphism $m_2$ and so is
itself a regular monomorphism. Since cokernel pairs of regular monomorphisms are pushouts along a regular monomorphism, they are stable 
under pullback, and so $i_2,j_2\colon A_2\rightrightarrows X_2$
is the cokernel pair of $m'_1$, and $e_2\colon X_2\to A_2$ is
the codiagonal.

Since $\ell$ is a regular monomorphism, we can form the pushout 
$$\xymatrix{
X_2 \ar[r]^{e_2} \ar[d]_{\ell} & A_2 \ar[d]^{k} \\ X_1 \ar[r]_{q} & Y.}$$
We shall see that the maps $qi,qj\colon X\to Y$ are the cokernel
pair of the union $m\colon A\to X$.

To do this, suppose that $u,v\colon X\to Z$ are given with 
$um=vm$, or equivalently with $um_1=vm_1$ and $um_2=vm_2$. 
Since $um_1=vm_1$, there is a unique
$w\colon X_1\to Z$ with $wi=u$ and $wj=v$. On the other hand
$w\ell i_2=wim_2=um_2=vm_2=wjm_2=w\ell j_2$ and so 
the morphism $w\ell$ out of the cokernel pair $X_2$ of $m'_1$ agrees
on the coprojections $i_2$ and $j_2$ of the cokernel pair, and therefore
factorizes through the codiagonal $e_2$, say as $w\ell=w'e_2$.
By the universal 
property of the pushout $Y$, there is a unique $w''\colon Y\to Z$
with $w''q=w$ and $w''k=w'$. Now $w''qi=wi=u$ and $w''qj=wj=v$,
and it is easy to see that $w''$ is unique with the property, thus
proving that $qi$ and $qj$ give the cokernel pair of $m$.

Since $qi$ and $qj$ are a cokernel pair, they have a common retraction
(the codiagonal), and so we can form their equalizer $n\colon B\to X$
as their intersection; this is of course a regular monomorphism.  
Since $qim=qjm$, there is a unique map $e\colon A\to B$ with $ne=m$;
we must show that $e$ is a stable epimorphism. In fact it will suffice to show that it is an epimorphism, since all these constructions are stable 
under pullback along a map into $X$, thus if $e$ is an epimorphism 
so will be any pullback. 

To see that $e$ is an epimorphism, first note that $m_1$ and $m_2$
factorize as $m_i=mn_i=nen_i$, and now $en_1$ and $en_2$ are regular
monomorphisms with union $e$. So once again we can factorize $e$ as $n'e'$,
where $n'$ is a regular monomorphism constructed as before. 

Observe that 
$e$ is an epimorphism if and only if its cokernel pair is trivial, which in 
turn is equivalent to the equalizer $n'$ of this cokernel pair being invertible. In an \rmadh category the regular monomorphisms are precisely the adhesive morphisms, and so are closed under composition. Thus the
composite $nn'$ of the regular monomorphisms $n$ and $n'$ is itself
a regular monomorphism. 
We conclude, using a standard argument, by showing that $n'$ is invertible and so that $e$ is an epimorphism. Observe that $qi$ and $qj$ are 
the cokernel pair of 
$m=ne=nn'e'$, as well as of their equalizer $n$. But then they
are also certainly the cokernel pair of $nn'$. Since this last is a regular
monomorphism, it must be the equalizer of $qi$ and $qj$, but then it
must be (isomorphic to) $n$, so that $n'$ is invertible as claimed.
\endproof
 
\begin{proposition}\label{prop:stably-jointly-epi}
If $m\colon C\to A$ is a regular monomorphism and $f\colon C\to B$ an
arbitrary morphism in an \qadh category, then the 
induced maps $\delta\colon A\to A_2$ and $m_2\colon C_2\to A_2$ are 
stably jointly epimorphic.  
\end{proposition}

\proof
Form, as in the proof of Lemma~\ref{lemma:basic}, the diagram 
$$\xymatrix{
C \ar[r]^{\gamma} \ar[d]_m & C_2 \ar[r]^{f_1} \ar[d]_{j} & C \ar[d]^{m} \\
A \ar[r]^{i} \ar[dr]_{\delta} & E \ar[r]^{h_1} \ar[d]^{k} & A \ar[d]^{1} \\
& A_2 \ar[r]_{g_1} & A }$$
with all squares pushouts, and note that $k\colon E\to A_2$ is the union 
of the regular monos $\delta\colon A\to A_2$ and $m_2\colon C_2\to A_2$. Thus
we can factorize it as a stable epimorphism $e\colon E\to\overline{E}$ followed
by a regular monomorphism $\overline{k}\colon \overline{E}\to A_2$.
Since $e$ is an epimorphism, the square
$$\xymatrix{
\overline{E} \ar[r]^{g_1\overline{k}} \ar[d]_{\overline{k}} & A \ar[d]^{1} \\
A_2 \ar[r]_{g_1} & A }$$
is also a pushout. Since $\overline{k}$ is a regular monomorphism,
the square is also a pullback, but then $\overline{k}$ is a pullback
of an invertible morphism (an identity), so is invertible.  Thus $k$ 
is a stable epimorphism, and so $\delta$ and $m_2$ are stably jointly 
epimorphic as claimed. 
\endproof

Let \C be a small \rmadh category. 
We have already described a topology  $j$
for which the representables are $j$-sheaves, and any $j$-sheaf
$F\colon\C\op\to\Set$  will 
send each pushout  along a regular monomorphism to a pullback, provided
that the induced morphsims $Fm_2$ and $F\delta$ are jointly monic. 
We therefore seek a topology $k$ containing $j$, with the property that
each representable is $k$-separated and each pair $\{m_2,\delta\}$ is 
a covering family. 

There are various possible choices for such a $k$, but perhaps it is
easiest to take  the topology whose covering families are 
generated by those in $j$ along with all pullbacks of those of the 
form $\{m_2,\delta\}$.  Since each such pair $\{m_2,\delta\}$ 
is stably jointly epimorphic by Proposition~\ref{prop:stably-jointly-epi},
the representables will indeed be $k$-separated. On the other hand, 
each $k$-separated
$F$ certainly has $Fm_2$ and $F\delta$ jointly monomorphic, and so
the $j$-sheaf condition for such an $F$ becomes the condition that $F$ send
pushouts along regular monomorphisms to pullbacks.

This gives:

\begin{theorem}\label{thmD}
Any small \qadh category \C has a full embedding,
preserving pushouts along regular monomorphisms and all existing
limits, into a Grothendieck quasitopos. 
\end{theorem}

Conversely, it is well-known that every quasitopos satisfies the 
conditions in the definition of \qadh category, hence so will any
full subcategory closed under the relevant limits and colimits. Thus
we have:

\enlargethispage{\baselineskip}

\begin{corollary}\label{corD}
Let \C be a small category with all pullbacks and with all pushouts along 
regular monomorphisms. Then \C is \qadh just when it has a full 
embedding, preserving the given structure, into a quasitopos.
\end{corollary}

 \bibliographystyle{plain}

\end{document}